\documentclass[11pt]{article}

\usepackage[centertags]{amsmath}
\usepackage[dvips]{graphicx}
\usepackage{amsfonts,amssymb,amsthm,enumerate,newlfont}
\usepackage{latexsym}
\usepackage{color}

\newcommand{\dem}{\noindent{\bf Proof \quad}}

\newtheorem {teo} {\bf Theorem\,} [section]
\newtheorem {prop} [teo] {\bf Proposition}
\newtheorem {coro} [teo]{\bf Corollary}
\newtheorem {lema} [teo] {\bf Lemma}
\newtheorem {defn} [teo] {\bf Definition}

\newtheorem {Remark} [teo] {\bf Remark}

\newcommand{\be}{\begin{eqnarray}}
\newcommand{\ee}{\end{eqnarray}}
\newcommand{\benn}{\begin{eqnarray*}}
\newcommand{\eenn}{\end{eqnarray*}}
\newcommand{\bse}{\begin{equation}}
\newcommand{\ese}{\end{equation}}
\newcommand{\bsenn}{\begin{displaymath}}
\newcommand{\esenn}{\end{displaymath}}
\numberwithin{equation}{section}

\title{Lower semicontinuity of global attractors for a class of evolution equations type neural fields in a bounded
domain}

\begin{document}

\author{\normalsize  Severino Hor\'acio da
Silva $^{1}$
\thanks{Partially supported by CAPES/CNPq-Brazil grant
Casadinho/Procad 552.464/2011-2 and INCTMat 5733523/2008-8.}
\hspace{0.1cm} \\{\scriptsize ${}^{1}$ Unidade Acad\^emica de
Matem\'atica e Estat\'istica UAME/CCT/UFCG}
\\{\scriptsize Rua Apr\'igio Veloso, 882,  Bairro Universit\'ario CEP 58429-900, Campina Grande-PB,
Brasil.}
\\{\scriptsize E-mail: horaciosp@gmail.com}
} \maketitle
\section*{Abstract}
In this work we consider the nonlocal evolution equation
$$
\frac{\partial u(w,t)}{\partial t}=-u(w,t)+
\int_{S^{1}}J(wz^{-1})f(u(z,t))dz+ h, \,\,\, h  > 0
$$
which arises in models of neuronal activity, in $L^{2}(S^{1})$,
where $S^{1}$ denotes the unit sphere. We obtain stronger results on
existence of global attractors and Lypaunov functional than the
already existing in the literature. Furthermore, we prove the
result, not yet known in the literature, of lower semicontinuity of
global attractors with respect to connectivity function $J$.
\noindent{ \it{Keywords:} Neural fields; Lypaunov functional; Lower
semicontinuity of attractors}

\section{Introduction} \label{intro}

We consider initially the nonlocal evolution equation proposed by
Wilson and Cowan in \cite{Wilson}, which is used to model neuronal
activity, that is,
\begin{equation}
\frac{\partial v(x,t)}{\partial t}=-v(x,t)+
\int_{\mathbb{R}}\widetilde{J}(x-y)f(v(y,t))dy+ h, \,\,\, h  > 0.
\label{WC}
\end{equation}
In (\ref{WC}), $v(x,t)$ is a real function on $\mathbb{R}\times
\mathbb{R}_{+},$ $\widetilde{J} \in C^{1}(\mathbb{R})$ is a non
negative even function supported in the interval $[-1,1]$, $f$ is a
non negative nondecreasing function and $h$ is a positive constant.

In this model, $v(x,t)$ denotes the mean membrane potential of a
patch of tissue located at position $x$ at time $t\geq 0$. The
connection function $\widetilde{J}$ determines the coupling between
the elements at position $x$ with the element at position $y$. The
non negative nondecreasing function $f(v)$ gives the neural firing
rate, or averages rate at which spikes are generated, corresponding
to an activity level $v$. The parameter $h$ denotes a constant
external stimulus applied uniformly to the entire neural field. We
say that the neurons at point $x$ is active if $S(x,t)>0$, where
$S(x,t)=f(v(x,t))$ is the firing rate of a neuron at position $x$ at
time $t$.

Proceeding as in \cite{Silva3}, it is easy to see that the Cauchy
problem for (\ref{WC}) is well posed in the space of continuous
bounded functions, $C_{b}(\mathbb{R})$, and that the subspace
$\mathbb{P}_{2\tau}$ of $2\tau$-periodic functions is invariant.
Thus, defining $\varphi:\mathbb{R} \rightarrow S^{1}$ by
$$
\varphi(x)=exp\left(i\frac{\pi}{\tau}x\right)
$$
and, for a $2\tau$ periodic function, $v$, defining
$u:S^{1}\rightarrow \mathbb{R}$ by $ u(\varphi(x))=v(x), $ and, in
particular, writing
$$
J(\varphi(x))=\widetilde{J}^{\tau}(x),
$$
where $\widetilde{J}^{\tau}$  denotes the $2 \tau$ periodic
extension of the restriction of $J$ to interval $[-\tau, \tau],$ for
some $\tau >1$, we obtain that: a function $v(x,t)$ is a $2\tau$
periodic solution  of (\ref{WC}) if and only if
$u(w,t)=v(\varphi^{-1}(w),t)$ is a solution of the equation
(\ref{1.1}) below:
\begin{equation}
\frac{\partial u(w,t)}{\partial t}=-u(w,t)+ J*(f\circ u)(w,t)+ h,
\,\,\, h  > 0, \label{1.1}
\end{equation}
where the $*$ above denotes convolution product in $S^{1}$, that is,
$$
(J*m)(w)=\int_{S^{1}}J(wz^{-1})m(z)dz,
$$
with $dz=\frac{\tau}{\pi}d\theta$, where $d\theta$ denotes
integration with respect to arc length.

In the literature, there are already several works dedicated to the
analysis of this model (see, for example, \cite{Amari}, \cite{Chen},
\cite{Ermentrout}, \cite{Ermentrout2}, \cite{Kishimoto},
\cite{Krisner}, \cite{Kubota}, \cite{Laing}, \cite{Rubin},
\cite{Silva}, \cite{Silva2}, \cite{Silva3} and \cite{Silva4}). Most
of these works have concerned with the existence and stability of
characteristic solutions, such as localized excitation (see, for
example, \cite{Amari}, \cite{Kishimoto} and \cite{Rubin})  or
traveling front (see, for example, \cite{Chen}, \cite{Ermentrout} and
\cite{Ermentrout2}). Also there are already some works on the global
dynamics of this model, (see, for example, \cite{Kubota},
\cite{Silva}, \cite{Silva2}, \cite{Silva3} and \cite{Silva4}).
However, the proof of the lower semicontinuity of global attractors
is not yet known, and this proof cannot be given by conventional
methods, since we cannot assume that equilibria are all hyperbolic,
leaving this property far more attractive from the point of view of
mathematical difficulty.

For the sake of clarity and future reference, it is convenient to
start with the hypotheses below used in \cite{Silva3} and
\cite{Silva4}.

\par\noindent {\bf (H1)} The function $f\in C^{1}(\mathbb{R})$,
$f'$ locally Lipschitz and
\begin{equation}
0<f'(r)<k_{1}, \, \forall \,\, r\in \mathbb{R},\label{1.2}
\end{equation}
for some positive constant $k_{1}$.

\par\noindent {\bf (H2)} $f$ is a nondecreasing function taking values between $0$ and $S_{max}>0$ and  satisfies, for $0\leq s\leq S_{max}$,
$$
\left|\int_{0}^{s}f^{-1}(r)dr\right|< L <\infty.
$$

From (H1) follows that
\begin{equation}
|f(x)-f(y)|\leq k_{1}|x-y|, \,\,\, \forall \, x,y \in
\mathbb{R},\label{1.3}
\end{equation}
and, in particular, there exists constant $k_{2}\geq 0$ such that
\begin{equation}
|f(x)|\leq k_{1}|x|+k_{2}.\label{1.4}
\end{equation}

In \cite{Silva3} and \cite{Silva4}, to obtain results on global
attractors and Lyapunov functional, besides the hypotheses (H1) and
(H2) above, it is assumed the hypothesis $k_{1}\|J\|_{L^{1}}<1$.
Under this assumption, the map $\Psi:L^{2}(S^{1})\rightarrow
L^{2}(S^{1})$ given by
$$
\Psi(u):=J*(f\circ u)+h
$$
is a contraction. Hence equation $(\text{\bf {P}})_{J}$ bellow has
an unique equilibrium $\bar{u}$, which can leave the attractor to
the trivial case of only one point.

In this paper, we organize the results as follows. In Section 2, we
conclude that the hypothesis $k_{1}\|J\|_{L^{1}}<1$ is not required
to obtain the results from \cite{Silva3} and \cite{Silva4} on global
attractors and Lyapunov functional. Therefore, we obtain (see
Theorem \ref{attractor} and Proposition \ref{gradient}) stronger
results in this direction. In Section \ref{lowersemi}, using the
same techniques of \cite{Severino2}, we prove the property of lower
semicontinuity of the attractors. To the extent of our knowledge,
with the exception of \cite{Severino2},
  the proofs of this
  property available  in the literature assume that
  \emph{ the
 equilibrium points are all hyperbolic} and therefore isolated
 (see for
 example \cite{AC2}, \cite{Bezerra}, \cite{OPP} and \cite{PP}).
  However,  this property cannot  hold true  in our case, due to the symmetries
 present in the equation.
   In fact,  it is   a consequence of these symmetries that
  the nonconstant equilibria
  arise in families and,  therefore, cannot be
  hyperbolic. This increases the difficulty and the interest of the problem, since we cannot use
  results of the type Implicit Function Theorem to prove the continuity of
  equilibria. To overcome this difficulty we have to  replace the  hypothesis of
hyperbolicity by  \emph{normal hyperbolicity} of curves of
equilibria. We  then used results
 of  \cite{Bates} on the permanence of \emph{normally hyperbolic invariant
 manifolds} and use  one result of \cite{Silva5}  of continuity properties of the local unstable manifolds of the
 curves of equilibria. Finally, in Section 4, we illustrate our
 results with a concrete example, which satisfies all hypotheses
 (H1)-(H4). This does not occur in \cite{Severino2} because there is no proven that the example
 satisfies the property that imply in normal hyperbolicity.

\section{Some remarks on global attractor and Lyapunov functional}

As proved in \cite{Silva4},
 under the hypothesis (H1),  the map
 \begin{equation} \label{mapF}
  F(u,J) = -u+J*(f(u))+
h
 \end{equation}
 is continuously Frechet
 differentiable in  $L^2(S^1) $ and, therefore,
 the equation
 \[  \hspace{50mm}
\frac{\partial u}{\partial t}= F(u,J) = -u+J*(f(u)) +h
   \hspace{30mm} (\text{\bf {P}})_{J}  \]
   generates   a $C^1$ flow in  $L^2(S^1) $ given, by the variation of constant
formula, by
$$
u(w,t)=e^{-t}u(w,0)+\int_{0}^{t}e^{-(t-s)}[J*(f\circ u)(w,s)+h]ds.
$$
From now on we denote this flow by $T_{J}(t)$ to make explicit
dependence on the parameter $J $.


Under hypothesis (H1), we proved in our previous work \cite{Silva3}
  that the Cauchy problem for  $(\text{\bf {P}})_{J}$, in
   $L^{2}(S^{1})$, is well posed and, assuming hypothesis (H1) and that $k_{1}\|J\|_{L^{1}}<1$, we proved  the existence and upper semicontinity of the global compact attractor in the sense of \cite{Hale}. Recently, in \cite{Silva4}, assuming the hypotheses (H1), (H2) and that $k_{1}\|J\|_{L^{1}}<1$, we prove that the flow of  $(\text{\bf {P}})_{J}$ is of class  $C^1$ and that it is gradient, in the sense of \cite{Hale},
   with Lyapunov functional ${\cal F}:L^{2}(S^{1}) \rightarrow \mathbb{R}$ given by
\begin{equation}
{\cal
F}(u)=\int_{S^{1}}\left[-\frac{1}{2}S(w)\int_{S^{1}}J(wz^{-1})S(z)dz+\int_{0}^{S(w)}f^{-1}(r)dr
-h S(w)\right]dw,\label{L1}
\end{equation}
where $S(w)=f(u(w))$.

It follows from Lemma below that we can obtain stronger versions of
Theorems 8 of \cite{Silva3} and Proposition 4.6 of \cite{Silva4},
eliminating the hypothesis $k_{1}\|J\|_{L^{1}}<1$ which is stronger
used these previous works.

\begin{lema}
Assume that (H1) and (H2) hold. Let $R = 2
\tau\|J\|_{\infty}S_{max}+h.$ Then the ball with center at the
origin of $L^{2}(S^{1})$ and radius $R\sqrt{2\tau}$ is an absorbing
set for the flow generated by  $(\text{\bf
{P}})_{J}$.\label{absorbing}
\end{lema}
\dem Let $u(w, t)$ be the solution of $(\text{\bf {P}})_{J}$ with
initial condition $u(w,0)$, then
$$
u(w,t)=e^{-t}u_{0}(w)+\int_{0}^{t}e^{-(t-s)}[J*(f\circ u)(w,s)+h]ds.
$$
Using hypothesis (H2) it follows that
\begin{eqnarray*}
|u(w,t)|&\leq& e^{-t}|u_{0}(w)|+\int_{0}^{t}e^{-(t-s)}|J*(f\circ u)(w,s)+h|ds\\
&\leq&e^{-t}|u_{0}(w)|+\int_{0}^{t}e^{-(t-s)}[2 \tau\|J\|_{\infty}S_{max}+h]\\
&\leq&e^{-t}|u_{0}(w)|+2 \tau\|J\|_{\infty}S_{max}+h\\
&=&e^{-t}|u_{0}(w)|+R.
\end{eqnarray*}
Hence,
\begin{eqnarray*}
\|u(\cdot,t)\|_{L^{2}}&\leq& \|e^{-t}|u_{0}|+R\|_{L^{2}}\\
&\leq&e^{-t}\|u_{0}\|_{L^{2}}+R\sqrt{2\tau}.
\end{eqnarray*}
Therefore, $u(\cdot,t)\in B(0,R+\varepsilon)$ for
$t>\ln\left(\frac{\|u_{0}\|_{L^{2}}}{\varepsilon}\right)$, and the
result is proved.\qed

From Lemma \ref{absorbing}, the Theorem 8 of \cite{Silva3} can be
rewritten as:
\begin{teo}
Suppose that the hypotheses (H1) and (H2) hold.  Then there exists a
global attractor ${\cal A}_{J}$ for the flow $T_{J}(t)$ in
$L^{2}(S^{1})$, which is contained in the ball of radius $(2
\tau\|J\|_{\infty}S_{max}+h)\sqrt{2\tau}$. \label{attractor}
\end{teo}

And from Theorem \ref{attractor}, the Proposition 4.6 of
\cite{Silva4}  can be rewritten as:
\begin{prop}
Assume that the hypothesis (H1) and (H2) hold. Then the flow
generated by equation $(\text{\bf {P}})_{J}$ is gradient, with
Lyapunov functional given by (\ref{L1}).\label{gradient}
\end{prop}

\section{Lower semicontinuity of the attractors} \label{lowersemi}

 As mentioned in the introduction, a additional difficulty we encounter in the proof of
lower semicontinuity is that,    due to the symmetries present in
our model,
 the  nonconstant equilibria are not isolated. In fact, as we will see
  shortly, the equivariance property of the map $F$ defined in  (\ref{mapF})
 implies that the nonconstant equilibria appear  in curves,
 (see Lemma \ref{Lema isolated}) and, therefore, cannot be hyperbolic
 preventing the use of tools like the Implicit Function Theorem to obtain
 their continuity with respect to parameters.

In this section we prove the lower semicontinuity property of
attractors, $\{{\cal A}_{J}\}$ at $J_{0}\in {\cal J},$ where
$$
{\cal J}=\{J\in C^{1}(\mathbb{R}), \, \mbox{even non negative,
supported in} \, [-1,1],\|J\|_{L^{1}}=1 \}.
$$

Let us recall that a family of subsets $\{{\cal A}_{J}\}$, is lower
semicontinuous at $J_{0}$ if
$$
dist( A_{J_{0}},A_{J})\longrightarrow 0, \,\,\mbox{as}\,\, J
\rightarrow J_{0}.
$$
where
\begin{equation} \label{defdist}
dist(A_{J_{0}}, A_{J})=\sup_{x\in A_{J_{0}}}dist(x,
A_{J})=\sup_{x\in A_{J_{0}}} \inf_{y\in A_{J}}\|x-y\|_{L^{2}}.
 \end{equation}

In order to obtain  the lower semicontinuity we will need  the
following  additional
   hypotheses:
\par\noindent
{\bf (H3)} For each $J_{0}\in {\cal J}$, the set $E$, of the
equilibria of $T_{J_{0}}(t)$, is such that $E=E_{1}\cup E_{2}$,
where
\par\noindent
{\bf (a)} the equilibria in $E_{1}$ are (constant)  hyperbolic
equilibria;
\par\noindent
{\bf (b)} the equilibria in  $E_{2}$ are  nonconstant and,
  for each $u_{0}\in E_{2}$, zero is simple eigenvalue of the
 derivative of $F$, with respect to $u$,  $D F_{u}(u_{0},J_{0}):L^{2}(S^{1})\rightarrow
L^{2}(S^{1})$, given by
$$
D F_{u}(u_{0},J_{0})v=-v+ J_{0}*(f'(u_{0})v).
$$
\par\noindent
 {\bf (H4)}  The function $f\in C^{2}(\mathbb{R})$.

  We start with some remarks on the spectrum of the linearization
 around equilibria.

 \begin{Remark}\label{autovalor}
 A simple computation  shows that, if $u_0$ is a nonconstant equilibria
  of $T_{J_{0}}(t)$ then   zero is always an eigenvalue of the operator
$$
DF_{u}(u_{0},J_{0})v=-v+ J_{0}*(f'(u_{0})v)
$$
with eigenfunction $u'_{0}$. Therefore, the hypothesis (H3)-b says
that
 we are in the `simplest'  possible situation for the linearization around nonconstant
 equilibria.
\end{Remark}

\begin{Remark} Let  $u_{0}\in E_{2}$.  It is easy to show that
 $D F_{u}(u_{0},J_{0})$ is a self-adjoint
operator with respect to  the inner product
$$
(u,v)=\int_{S^{1}}u(w)v(w)d(w).
$$
\label{Remark 3.1} Since
$$
v\rightarrow J_{0}*(f'(u_{0})v)
$$
is a  compact operator in  $L^{2}(S^{1})$, it  follows  from (H3)
that
 $$ \sigma(DF_{u}(u_{0},J_{0}))\backslash\{0\} $$
 contains only  real eigenvalues of finite multiplicity with
 $-1$ as  the unique possible accumulation point.
 \end{Remark}

Now we enunciate a result on the structure of the sets  of
 nonconstant equilibria. The proof of this result is omitted because it is very similar to the proof of the Lemma 3.1
 of \cite{Severino2}

\begin{lema} \label{curveseq}
Suppose that, for some $J_{0}\in {\cal J}$, (H1),  (H3) and (H4)
hold. Given $u\in E_{2}$ and $\alpha \in S^1$,
 define $\gamma(\alpha;u) \in L^{2}(S^1)$ by
$$
\gamma(\alpha;u)(w)=u(\alpha w), \,\, w \in S^{1}.
$$
Then $\Gamma=\gamma(S^{1};u)$ is a  closed, simple $C^{2}$
   curve  of equilibria of
  $T_{J_{0}}(t)$  which is  isolated in
  the set of equilibria, that is, no  point
of $\Gamma$ is an  accumulation point of $E_{J_{0}}\setminus
\Gamma$.\label{Lema isolated}
\end{lema}

\begin{coro}
Let $M$ a closed connected curve of equilibria in  $E_{2}$ and
$u_{0}\in M$. Then $M = \Gamma$, where
$\Gamma=\gamma(S^{1},u_{0})$.\label{Coro caract. das curvas}
\end{coro}
\dem Suppose that $ \Gamma \not\subset M $. Then there exist
equilibria in $M\setminus \Gamma$  accumulating at $u_{0}$
contradicting Lemma \ref{Lema isolated}. Therefore $ \Gamma
\subseteq M  $. Since $\Gamma$ is a simple closed curve, it follows
that $M=\Gamma$.\qed

   In order  to  prove our main result, we need
 some preliminary results, which we present  in the next three subsections.

\subsection{Continuity of the equilibria}

The upper semicontinuity of the equilibria is a consequence of the
upper semicontinuity of global attractors (see Theorem 11 of
\cite{Silva3}). The lower semicontinuity of the {\em hyperbolic}
equilibria is
 usually  obtained
 via the  Implicit Function Theorem.  However, this approach
 fails here since the
  equilibria may appear in families as we have shown in Lemma
 \ref{curveseq}.  To  overcome this difficulty, we  need
   the
concept of normal hyperbolicity, (see \cite{Bates}).

Recall that, if
 $T(t): X \to X $ is a semigroup,
 a set $ M \subset X$ is {\em invariant} under  $T(t)$ if
  $T(t) M = M$, for any $t > 0$.

\begin{defn}
Suppose that $T(t)$ is  a $C^1$  semigroup in a  Banach space $X$
and $M \subset X$ is an invariant manifold for  $T(t)$. We say that
$M$ is normally hyperbolic under $T(t)$ if
\par\noindent
{\bf (i)} for each $m\in M$ there is a decomposition
$$
X=X_{m}^{c}\oplus X_{m}^{u}\oplus X_{m}^{s}
$$
 by  closed subspaces with $X_{m}^{c}$ being  the tangent space to $M$ at
$m$.
\par\noindent
{\bf (ii)} for each $m\in M$ and $t\geq 0$, if $m_{1}=T(t)(m)$
$$
DT(t)(m)|_{X_{m}^{\alpha}}: X_{m}^{\alpha}\rightarrow
X_{m_{1}}^{\alpha}, \,\, \alpha = c, u, s
$$
and $DT(t)(m)|_{X_{m}^{u}}$ is an isomorphism from $X_{m}^{u}$ onto
$X_{m_{1}}^{u}$.
\par\noindent
{\bf (iii)} there is $t_{0}\geq 0$ and $\mu <1$ such that for all
$t\geq t_{0}$
\begin{equation}
\mu \inf \left\{\|DT(t)(m)x^{u}\| \, : \, x^{u}\in X_{m}^{u}, \,
\|x^{u}\|=1\} > \max\{1,
\|DT(t)(m)|_{X_{m}^{c}}\|\right\},\label{7.17}
\end{equation}
\begin{equation}
\mu \min \left\{1, \inf\{\|DT(t)(m)x^{c}\| \, : \, x^{c}\in
X_{m}^{c}, \, \|x^{c}\|=1\} \right\} >
\|DT(t)(m)|_{X_{m}^{s}}\|.\label{7.18}
\end{equation}\label{Def 7.2}
\end{defn}
The condition (\ref{7.17}) suggests that near $m\in M$, $T(t)$ is
expansive in the direction of $X_{m}^{u}$ and at rate greater than
on $M$, while (\ref{7.18}) suggests that $T(t)$  is contractive in
the direction of $X_{m}^{s}$, and at a rate greater than that on
$M$.

 The following result has been proved in \cite{Bates}.

\begin{teo} ({\bf Normal Hyperbolicity})
Suppose that $T(t)$ is a $C^{1}$ semigroup on a Banach space $X$ and
$M$ is a $C^{2}$ compact connected invariant manifold which is
normally hyperbolic under $T(t)$, ( that is (i) and (ii) hold and
there exists $0\leq t_{0}< \infty$ such that (iii) holds for all
$t\geq t_{0}$). Let $\widetilde{T}(t)$ be a $C^{1}$ semigroup on $X$
and $t_{1}>t_{0}$. Consider $N(\varepsilon)$, the
$\varepsilon$-neighborhood of $M$, given by
$$
N(\varepsilon)=\{m+x^{u}+x^{s}, \,\, x^{u}\in X^{u}_{m}, \,\,
x^{s}\in X^{s}_{m}, \,\, \|x^{u}\|, \, \|x^{s}\|< \varepsilon\}.
$$
Then, there exists $\varepsilon^{*}>0 $ such that for each
$\varepsilon < \varepsilon^{*}$, there exists $\sigma >0$ such that
if
$$
\sup_{u\in
N(\varepsilon)}\left\{\|\widetilde{T}(t_{1})u-T(t_{1})u\|+
\|D\widetilde{T}(t_{1})(u)-DT(t_{1})(u)\|\right\}< \sigma
$$
and
$$
\sup_{u\in N(\varepsilon)}\|\widetilde{T}(t)u-T(t)u\|< \sigma,
\mbox{for} \,\, 0\leq t\leq t_{1},
$$
 there is an unique compact connected invariant manifold of
class $C^{1}$, $\widetilde{M}$, in $N(\varepsilon)$. Furthermore,
$\widetilde{M}$ is normally hyperbolic under $\widetilde{T}(t)$ and,
for each $t\geq 0$, $\widetilde{T}(t)$ is a ${C^{1}}$-diffeomorphism
from  $\widetilde{M}$ to $\widetilde{M}$.
 \label{Teorema 7.3}
\end{teo}

\begin{prop}
Assume that the hypotheses (H1), (H2) and (H3) hold. Then, for each
$J \in {\cal J}$, any  curve of equilibria of $T_{J}(t)$ is a
normally hyperbolic manifold  under $T_{J}(t)$.\label{Prop 7.4}
\end{prop}
\dem Here we follow closely a proof of  \cite{Severino2}. Let $M$ be
a curve of equilibria of $T_{J}(t)$ and $m\in M$. From (H3) it
follows that
$$
Ker (DF_{u}(m,J))=span\{m'\}.
$$
Let $Y=\mathcal{R}(DF_{u}(m,J))$  the range of $DF_{u}(m,J)$. Since
$DF_{u}(m,J)$ is self-adjoint and Fredholm of index zero, it follows
from (H3) that
$$
\sigma(DF_{u}(u_{0},J)|_{Y})=\sigma_{u}\cup \sigma_{s},
$$
where $\sigma_{u}$, $\sigma_{s}$ correspond  to the positive and
 negative eigenvalues  respectively.

From (H1) and (H2), it  follows that $T_{J}(t)$ is a $C^{1}$
semigroup.  Consider the linear autonomous equation
\begin{equation} \label{DT}
\dot{v}=(DF_{u}(m,J)|_{Y})v.
 \end{equation}
  Then   $DT_{J}(t) v_0$   is the solution   of (\ref{DT})
  with initial condition $v_0$, that is  $DT_{J}(t)(m)v_0 =e^{(DF_{u}(m,J))t} v_0$. In
particular $DT_{J}(t)(m)|_{Y}\equiv
D(T_{J}(t)|_{Y})(m)=e^{(DF_{u}(m,J)|_{Y})t}$.

Let $P_{u}$ and $P_{s}$ be the spectral projections corresponding to
$\sigma_{u}$ and $\sigma_{s}$. The subspaces $X_{m}^{u}=P_{u}Y$,
$X_{m}^{s}=P_{s}Y$ are then invariant under
 $ D T_{J}(t) $ and the following estimates hold
  (see \cite{Daleckii},
p. 73, 81
).  \\
\begin{equation}
\|DT_{J}(t)|_{Y}v\| \leq Ne^{-\nu t}\|v\|, \,\,\mbox{for}\,\,\, v\in
X_{m}^{s} \,\,\mbox{and}\,\, t\geq 0, \label{Atract.exp}
\end{equation}
\begin{equation}
\|DT_{J}(t)|_{Y}v\| \leq Ne^{\nu t}\|v\|, \,\,\mbox{for}\,\,\, v\in
X_{m}^{u} \,\,\mbox{and}\,\, t\leq 0, \label{Rep.exp}
\end{equation}
for some positive constant $\nu$ and some constant $N>1$.

 It is clear that $DT_{J}(t) \equiv 0$ when
  restricted to
$X_{m}^{c}=span\{m'\}$.
  Therefore, we  have the decomposition
$$
L^{2}(S^{1})=X_{m}^{c}\oplus X_{m}^{u}\oplus X_{m}^{s}.
$$
 Since  $DF_{u}(m,J)|_{Y}$ is an isomorphism
$$
DF_{u}(m,J)|_{X_{m}^{\alpha}}:X_{m}^{\alpha}\rightarrow
X_{m}^{\alpha}, \,\, \alpha=u,s,
$$
is an isomorphism. Consequently,  the linear flow
$$
DT_{J}(t)(m)|_{X_{m}^{u}}:X_{m}^{u}\rightarrow X_{m}^{u}
$$
is also an isomorphism.

Finally, the estimates (\ref{7.17}) and (\ref{7.18}) follow from
 estimates (\ref{Atract.exp}) and (\ref{Rep.exp})
above. \qed

\begin{Remark}
For $u,v\in L^{2}(S^{1})$, from (\ref{1.2}) follows  that
\begin{eqnarray}
\|f'(u)v\|\leq k_{1}\|v\|_{L^{2}}\label{BoudendL2}.
\end{eqnarray}
\end{Remark}

\begin{prop}
Suppose that the hypotheses (H1)-(H2) hold. Let $DT_{J}(t)(u)$ be
the linear flow generated by the equation
$$
\frac{\partial v}{\partial t}=-v+ J_{0}*(f'(u_{0})v).
$$
Then, for a fixed $J_{0}\in {\cal J}$, we have
$$
\|T_{J}(t)u-T_{J_{0}}(t)u\|_{L^{2}(S^{1})}+ \|D T_{J}(t)(u) - D
T_{J_{0}}(t)(u)\|_{\mathcal{L}(L^{2}(S^{1}), \, L^{2}(S^{1}))}
\rightarrow 0 \,\, \mbox{as} \,\, \|J- J_{0}\|_{L^{1}} \rightarrow
0,
$$
 uniformly for $u$ in  bounded sets of $L^{2}(S^{1})$ and $t\in [0,b]$,
 $b<\infty$.
 \label{Prop 7.5}
\end{prop}
\dem From Lemma 10 of \cite{Silva3} it  follows that
$$
\|T_{J}(t)u-T_{J_{0}}(t)u\|_{L^{2}(S^{1})}\rightarrow 0 \,\,
\mbox{as} \,\, \|J- J_{0}\|_{L^{1}} \rightarrow 0,
$$
for $u$ in  bounded sets of $L^{2}(S^{1})$ and $t\in [0,b]$.

 By the variation of
constants formula, we have
$$
DT_{J}(t)(u)v=e^{-t}v+\int_{0}^{t}e^{-(t-s)}J*(f'(u)v)ds.
$$
Thus, using Young's inequality, we obtain
\begin{eqnarray*}
 \|D
T_{J}(t)(u)v - D T_{J_{0}}(t)(u)v\|_{L^{2}} &\leq&
\int_{0}^{t}e^{-(t-s)}\big\|(J-J_{0})*(f'(u)v) \big\|_{L^{2}}ds\\
&\leq& \int_{0}^{t}e^{-(t-s)}\|J-J_{0}\|_{L^{1}}\|f'(u)v\|_{L^{2}}.
\end{eqnarray*}
Using (\ref{BoudendL2}), it follows that
\begin{eqnarray*}
\|D T_{J}(t)(u)v - D T_{J_{0}}(t)(u)v\|_{L^{2}} &\leq&
k_{1}\|J-J_{0}\|_{L^{1}}\|v\|_{L^{2}}.
\end{eqnarray*}

Therefore
\begin{eqnarray*}
\|D T_{J}(t)(u) - D T_{J_{0}}(t)(u)\|_{\mathcal{L}(L^{2}(S^{1}),
\,\,L^{2}(S^{1}))} &=& \sup_{\|v\|=1}\|D T_{J}(t)(u)v - D
T_{J_{0}}(t)(u)v\|_{L^{2}(S^{1})} \\
&\leq&
\sup_{\|v\|=1}k_{1}\|J-J_{0}\|_{L^{1}}\|v\|_{L^{2}}\\
&=&C(J),
\end{eqnarray*}
with $C(J)\rightarrow 0$, as $\|J-J_{0}\|_{L^{1}}\rightarrow 0$.
This completes the proof.\qed

The proof of the theorem below follows closely the proof of
Theorem 3.4 of \cite{Severino2}.

\begin{teo}
Suppose that the hypotheses (H1)-(H4) hold. Then the set $E_{J}$ of
the equilibria of $T_{J}(t)$ is lower semi-continuous with respect
to $J$ at $J_{0}$.\label{Teorema 7.6}
\end{teo}
\dem  The continuity of the  constant equilibria  follows from the
Implicit Function Theorem and the hypothesis of hyperbolicity.

 Suppose now that  $m$ is  a  nonconstant equilibrium  and let
 $\Gamma= \gamma(\alpha;m) $  be the isolated curve
 of equilibria containing $m$ given by Lemma \ref{Lema isolated}.
  We want to show that, for every
$\varepsilon
>0$, there exists $\delta>0$  so that,  if $J\in {\cal J}$ , there exists
$\Gamma_{J}\in E_{J}$ such that $\Gamma \subset
\Gamma_{J}^{\varepsilon}$, where $\Gamma_{J}^{\varepsilon}$ is the
$\varepsilon$-neighborhood of $\Gamma_{J}$.

From Lemma \ref{Lema isolated} and Propositions \ref{Prop 7.4} and
\ref{Prop 7.5},  the assumptions of the Normal Hyperbolicity Theorem
are met. Thus,  given $\varepsilon >0$, there is $\delta
>0$ such that, if $\|J-J_{0}\|_{L^{1}}<\delta$ there is an unique $C^{1}$ compact connected invariant
manifold  $\Gamma_{J}$ normally hyperbolic under $T_{J}(t)$, such
that $\Gamma_{J}$ is $\varepsilon$-close
 and $C^{1}$-diffeomorphic  to $\Gamma$.

Since $T_{J}(t)$ is gradient and $\Gamma_{J}$ is compact, there
exists at least one equilibrium $m_{J} \in \Gamma_{J}$.
 In fact, the $\omega$ limit of any $u \in \Gamma_{J}$ is nonempty and
  belongs to $\Gamma_{J}$ by invariance.
  From Lemma 3.8.2 of
\cite{Hale},  it must contain an equilibrium. Since $\Gamma_{J}$ is
$\varepsilon$-close to $\Gamma$, there exists $m\in \Gamma$ such
that $ \|m-m_{J}\|_{L^2(S^1)}<\varepsilon. $

 Let $\tilde{\Gamma}_{J}$ be the curve of equilibria given by
  $\tilde{\Gamma}_{J}\equiv\{\gamma(\alpha;m_{J}), \,
 \alpha \in S^{1}\}$ which is
 a  normally hyperbolic invariant manifold under $T_{J}(t)$
  by Proposition \ref{Prop 7.4}.
Then, for each $\alpha \in S^{1}$, we have
\begin{eqnarray*}
\|\gamma(\alpha;m_{J})-\gamma(\alpha;m)\|_{L^{2}}^{2}
&=&\int_{S^{1}}|\gamma(\alpha;m_{J})(w)-\gamma(\alpha;m)(w)|^{2}dw\\
&=&\int_{S^{1}}|m_{J}(\alpha w)-m(\alpha w)|^{2}dw \\
& = & \|m_{J}-m\|_{L^{2}}.\\
\end{eqnarray*}
Thus
\begin{eqnarray*}
\|\gamma(\alpha;m_{J})-\gamma(\alpha;m)\|_{L^{2}} &=&
\|m_{J}-m\|_{L^{2}}\\
&<& \varepsilon.
\end{eqnarray*}
and  $\Gamma$ is $\varepsilon$-close to $ \tilde{\Gamma_{J}}.$
 Since  there are only a finite number of curves of equilibria the result
 follows immediately. \qed\\

%
%
%
%
%
%
%
%
%

\subsection{Existence and continuity of the local unstable manifolds}

 Let us return to   equation $(\text{\bf {P}})_{J}$.
 Recall that the \emph{unstable set}
 $W_{J}^{u}= W_{J}^{u}(u_{J})$ of  an equilibrium  $u_{J}$
 is the set of initial conditions $\varphi$ of  $(\text{\bf {P}})_{J}$,
  such that $T_{J}(t) \varphi$ is defined for all $t \leq 0$ and
 $T_{J}(t) \varphi \to u_{J} $ as $t \to -\infty$.
  For a  given neighborhood $V$ of $u_{J}$,  the set
 $W_{J}^{u} \cap V$ is called a \emph{local unstable set} of
 $u_{J}$.

   Using results of \cite{Silva5}   we now show that the local unstable sets
  are actually Lipschitz manifolds   in a sufficiently small neighborhood
  and  vary continuously with
 $J$. More precisely, we have

  \begin{lema} \label{contunstball}
   If
   $ u_{0}  $ is  a fixed  equilibrium of  $(\text{\bf {P}})_{J}$
 for
  $J = J_0$, then
 there is a
  $\delta> 0$ such
   that,  if  $ \| J - J_0\|_{L^{1}} + \| u_0 - u_{J}\|_{L^{2}} < \delta$
   and
 \[  U_{J}^\delta :=  \{ u \in  W_{J}^{u}(u_{J}) \  : \
 ||u-u_{J}  ||_{L^2}< \delta \}
  \]
  then  $U_{J}^\delta$ is a Lipschitz manifold and

  \[ \textrm{dist}( U_{J}^\delta , U_{J_{0}}^\delta     ) +
dist( U_{J_{0}}^\delta, U_{J}^\delta   )
  \to 0  \quad \text{as}  \quad
  \| J - J_0\|_{L^{1}} + \| u_0 - u_{J}\|_{L^{2}}    \to  0, \]
  with  \textrm{dist}  defined as in (\ref{defdist}).
   \end{lema}

  \dem
   As already mentioned in the previous section,
 assuming  the hypothesis (H1),  the
map  $F:L^{2}(S^{1})\times {\cal J} \rightarrow L^{2}(S^{1})$,
 $$
 F(u, J) = -u + J*(f(u))+h
$$
  defined by the   right-hand side of $(\text{\bf {P}})_{J}$
  is continuously
 Frechet differentiable.
  Let $ u_{J}$ be an equilibrium of
 $(\text{\bf {P}})_{J}$.
   Writing $u = u_{J} + v$, it follows that $u$ is a solution
 of $(\text{\bf {P}})_{J}$ if and only if $v$ satisfies

 \begin{equation}
\frac{\partial v}{\partial t}=L(J )v + r(u_{J},v,J),\label{7.23}
\end{equation}
where $L(J)v= \frac{\partial}{\partial u} F( u_{J},J ) =
 -v+ J*(f'(u_{J})v)$ and
 $  r(u_{J},v,J) =   F( u_{J}+v,J ) -
   F( u_{J},J ) -
    L(J) v.   $
We rewrite  equation (\ref{7.23}) in the form
\begin{equation}
\frac{\partial v}{\partial t}=L(J_{0})v+ g(v, J),\label{7.26}
\end{equation}
where $ g(v,J) = [L(J) - L(J_{0})]v + r(u_{J},v, J) $ is  the ``non
linear part" of (\ref{7.26}).
 Note that now  the ``linear part'' of (\ref{7.26})
 does not  depend on
 the parameter $J$, as required by Theorems 2.5 and 3.1 from \cite{Silva5}.

Note that
\begin{eqnarray*}
\|[L(J)-L(J_{0})]v\|_{L^{2}} \leq \|(J-J_{0})*(f'(u_{J})v)\|_{L^{2}}
+ \|J_{0}*[f'(u_{J})-f'(u_{J_{0}})]v\|_{L^{2}}.
\end{eqnarray*}
But, using Holder inequality, we have
\begin{eqnarray*}
|J_{0}*[f'(u_{J})-f'(u_{J_{0}})](w)v(w)|&\leq& \int_{S^{1}}J_{0}(wz^{-1})|f'(u_{J}(z))-f'(u_{J_{0}}(z))||v(z)|dz\\
&\leq& \|J_{0}\|_{\infty}\int_{S^{1}}|f'(u_{J}(z))-f'(u_{J_{0}}(z))||v(z)|dz\\
&\leq& \|J_{0}\|_{\infty}
\|f'(u_{J})-f'(u_{J_{0}})\|_{L^{2}}\|v\|_{L^{2}}.
\end{eqnarray*}
Thus, remembering that we are assuming the notation of our previous
work (\cite{Silva3} and \cite{Silva4}), where the measure of $S^{1}$
is $2\tau$ , we obtain
\begin{eqnarray*}
\|J_{0}*[f'(u_{J})-f'(u_{J_{0}})]v\|_{L^{}2} \leq \sqrt{2\tau}
\|J_{0}\|_{\infty} \|f'(u_{J})-f'(u_{J_{0}})\|_{L^{2}}\|v\|_{L^{2}}.
\end{eqnarray*}
Hence, using  Young inequality and hypothesis (H1), we have
\begin{eqnarray*}
\|[L(J)-L(J_{0})]v\|_{L^{2}} \leq
k_{1}\|J-J_{0}\|_{L^{1}}\|v\|_{L^{2}} +\sqrt{2\tau}
\|J_{0}\|_{\infty}
\|f'(u_{J})-f'(u_{J_{0}})\|_{L^{2}}\|v\|_{L^{2}}.\label{II}
\end{eqnarray*}

But, keeping $u_{J_{0}}\in L^2(S^1)$, from Theorem \ref{Teorema
7.6},  follows that $\|u_{J} - u_{J_{0}}\|_{L^{2}}\rightarrow 0$, as
$\|J-J_{0}\|_{L^{1}}\rightarrow 0$. It follows that $u_{J}(w)
\rightarrow u_{J_{0}}(w)$ almost everywhere in $S^1$. From (H1)
follows that, there exists $M>0$
such that
$$
|f'(u_{J}(w))-f'(u_{J_{0}}(w))|\leq M|u_{J}(w)-u_{J_{0}}(w)|, \quad
\text{almost everywhere}.
$$
Then
\begin{eqnarray}
\|f'\circ u_{J} - f'\circ u_{J_{0}}\|_{L^2}^2&=& \int_{S^1}|f'(u_{J}(w))-f'(u_{J_{0}}(w))|^2dw\nonumber \\
&\leq& \int_{S^1}M^2|u_{J}(w)-u_{J_{0}}(w)|^2dw\nonumber \\
&=& M^2\|u_{J}-u_{J_{0}}\|_{L^2}^2.\label{III}
\end{eqnarray}
Therefore, using (\ref{III}), we obtain
\begin{eqnarray}
\|[L(J)-L(J_{0})]v\|_{L^{2}} &\leq&
k_{1}\|J-J_{0}\|_{L^{1}}\|v\|_{L^{2}} + \sqrt{2\tau}
\|J_{0}\|_{\infty} M\|u_{J}-u_{J_{0}}\|_{L^2}\|v\|_{L^{2}}
.\label{IV}
\end{eqnarray}
%
Now, note that,
\begin{eqnarray*}
r(u_{J},v,J)-r(u_{J_{0}},v,J_{0})&=&F(u_{J}+v,J)-F(u_{J},J)-L(J)v\\
&-&F(u_{J_{0}}+v,J_{0})+F(u_{J_{0}}, J_{0})+L(J_{0})v\\
&=&J*f(u_{J}+v)-J*f(u_{J})+J_{0}*f(u_{J_{0}})-J_{0}*f(u_{J_{0}}+v)\\
&-&[L(J)-L(J_{0})]v\\
&=&J*[f(u_{J}+v)-f(u_{J})]
+J_{0}*[f(u_{J_{0}})-f(u_{J_{0}}+v)]\\
&-&[L(J)-L(J_{0})]v.
\end{eqnarray*}
But
\begin{eqnarray*}
J*[f(u_{J}+v)-f(u_{J})]=J*[f'(\bar{v})v]
\end{eqnarray*}
and
\begin{eqnarray*}
J_{0}*[f(u_{J_{0}}+v)-f(u_{J_{0}})]=J*[f'(\bar{\bar{v}})v],
\end{eqnarray*}
for some $\bar{v}$ in the segment defined by $J*f(u_{J})$ and
$J*f(u_{J}+v)$ and for some $\bar{\bar{v}}$ in the segment defined
by $J_{0}*f(u_{J_{0}})$ and $J_{0}*f(u_{J_{0}}+v)$. Then
\begin{eqnarray*}
J*[f(u_{J}+v)-f(u_{J})]
+J_{0}*[f(u_{J_{0}})-f(u_{J_{0}}+v)]&=& J*[f'(\bar{v})v]-J_{0}*[f'(\bar{\bar{v}})v]\\
&=&J*[f'(\bar{v})]v-J_{0}*[f'(\bar{v})v]\\
&+&J_{0}*[f'(\bar{v})v]-J_{0}*[f'(\bar{\bar{v}})v]\\
&=&(J-J_{0})*f'(\bar{v})v\\
&+&J_{0}*[f'(\bar{v})-f'(\bar{\bar{v}})]v.
\end{eqnarray*}
Thus
\begin{eqnarray*}
r(u_{J},v,J)-r(u_{J_{0}},v,J_{0})&=&(J-J_{0})*f'(\bar{v})v
+J_{0}*[f'(\bar{v})-f'(\bar{\bar{v}})]v\\
&+&[L(J_{0})-L(J)]v.
\end{eqnarray*}
Hence
\begin{eqnarray*}
&\|&r(u_{J},v,J)-r(u_{J_{0}},v,J_{0})\|_{L^{2}}\\
&\leq& \|J-J_{0}\|_{L^{1}}\|f'(\bar{v})v\|_{L^{2}}
+\|J_{0}*[f'(\bar{v})-f'(\bar{\bar{v}})]v\|_{L^{2}}
+\|[L(J)-L(J_{0})]v\|_{L^{2}}.
\end{eqnarray*}
But, from hypothesis (H4), there exists $M>0$ such that
$$
|f'(\bar{v})(z)-f'(\bar{\bar{v}})(z)|\leq
M|\bar{v}(z)-\bar{\bar{v}}(z)|, \,\, \forall \,\, z\in S^{1},
$$
thus
\begin{eqnarray*}
|J_{0}*[f'(\bar{v})(w)-f'(\bar{\bar{v}})(w)]v(w)|
&\leq& \int_{S^{1}} J_{0}(wz^{-1})|f'(\bar{v})(z)-f'(\bar{\bar{v}})(z)||v(z)|dz\\
&\leq& \int_{S^{1}}\|J_{0}\|_{\infty}M|\bar{v}(z)-\bar{\bar{v}}(z)||v(z)|dz\\
&\leq&
\|J_{0}\|_{\infty}M\|\bar{v}-\bar{\bar{v}}\|_{L^{2}}\|v\|_{L^{2}}.
\end{eqnarray*}
Then, remembering again that we are assuming the notation of our
previous work (\cite{Silva3} and \cite{Silva4}), where the measure
of $S^{1}$ is $2\tau$, we obtain
\begin{eqnarray}
\|J_{0}*[f'(\bar{v})-f'(\bar{\bar{v}})]v\|_{L^{2}} \leq
\|J\|_{\infty}M\sqrt{2\tau}\|\bar{v}-\bar{\bar{v}}\|_{L^{2}}\|v\|_{L^{2}}
,\label{V}
\end{eqnarray}
Thus, using (\ref{BoudendL2}), (\ref{IV}) and (\ref{V}), and the
fact that $\|\bar{v}-\bar{\bar{v}}\|_{L^{2}}\rightarrow 0$, as
$\|J-J_{0}\|_{L^{1}}\rightarrow 0$, follows that
\begin{eqnarray}
\|r(u_{J},v,J)-r(u_{J_{0}},v,J_{0})\|_{L^{2}} &\leq&
k_{1}\|J-J_{0}\|_{L^{1}}\|v\|_{L^{2}}
+\|J\|_{\infty}M\sqrt{2\tau}\|\bar{v}-\bar{\bar{v}}\|_{L^{2}}\|v\|_{L^{2}} \nonumber\\
&+&k_{1}\|J-J_{0}\|_{L^{1}}\|v\|_{L^{2}} +
\sqrt{2\tau} \|J_{0}\|_{\infty} M\|u_{J}-u_{J_{0}}\|_{L^2}\|v\|_{L^{2}}\nonumber\\
&=&C_{1}(J)\|v\|_{L^{2}},\label{VI}
\end{eqnarray}
with $C_{1}(J)\rightarrow 0$, as $\|J-J_{0}\|_{L^{1}}\rightarrow 0$.

Now, since
$$
g(v,J)-g(v,J_{0})=[L(J)-L(J_{0})]v +r(u_{J},v,J)-r(u_{J_{0}},v,
J_{0}),
$$
using (\ref{IV}) and (\ref{VI}), we obtain
\begin{eqnarray}
\|g(v,J)-g(v,J_{0})\|&\leq& \|L(J)-L(J_{0}\|_{L^{2}}+\|r(u_{J},v,J)-r(u_{J_{0}},v,J_{0})\|_{L^{2}}\nonumber\\
&\leq&C_{2}(J), \label{VII}
\end{eqnarray}
where $C_{2}(J)\rightarrow 0$  as $\|J-J_{0}\|_{L^{2}} \rightarrow
0$.

 In a similar way, we obtain for any $ v_{1}, v_{2}$ with
 $ || v_{1}||_{L^2(S^1)} $ and  $ || v_{2}||_{L^2(S^1)} $ smaller than $\rho$




\begin{equation}
\|g(v_{1},J)-g(v_{2}, J)\|_{L^{2}}  \leq
\nu(\rho)\|v_{1}-v_{2}\|_{L^{2}} ,\label{estf(v,w)}
\end{equation}
where  $\nu(\rho) \to 0$, as $\rho \to 0$.

Therefore, the conditions of Theorems 2.5 and 3.1 from \cite{Silva5}
are satisfied and  we obtain the existence of locally
 invariant sets  for  (\ref{7.26}) near the origin,
  given as graphs of Lipschitz
 functions which depend continuously on the parameter $J$
 near $J_0$.  Using uniqueness of solutions, we can easily
 prove that these sets coincide with the
  local unstable manifolds
 of  (\ref{7.26}).

 Now, noting that the translation
$$
u\rightarrow (u-u_{J})
$$
sends an equilibrium $u_{J}$ of $(\text{\bf {P}})_{J}$
 into the origin
(which is an equilibrium of (\ref{7.26})),
 the results claimed follow immediately.
 \qed

 Using the compactness of the set of equilibria, one can obtain an
 `uniform version' of Lemma  \ref{contunstball} that will be
 needed later.

  \begin{lema} \label{contunstballunif}
    Let  $J = J_0$ be fixed.
 Then, there is a  $\delta> 0$ such
   that, for any equilibrium
   $ u_{0}  $  of   $(\text{\bf {P}})_{J_0}$,
    if   $ \| J - J_0\|_{L^{1}} + \| u_0 - u_{J}\|_{L^{2}} < \delta$
   and
 \[  U_{J}^\delta :=  \{ u \in  U_{J}(u_{J}) \  : \
 ||u-u_{J}  ||_{L^2(S^1)}< \delta \}
  \]
  then  $U_{J}^\delta$ is a Lipschitz manifold and
 \[ \sup_{u_0 \in E_{J_0}} \textrm{dist}( U_{J}^\delta , U_{J_{0}}^\delta     ) +
dist( U_{J_{0}}^\delta, U_{J}^\delta   )
  \to 0  \quad \text{as}  \quad
  \| J - J_0\|_{L^{1}} + \| u_0 - u_{J}\|_{L^{2}}    \to  0, \]
  with  \textrm{dist}  defined as in (\ref{defdist})
   \end{lema}
   \dem  From Lemma \ref{contunstball}, we know that, for any
   $ u_{0} \in E_{J_0} $, there is a $\delta = \delta(u_0)$ such that
   $ U_{J}^\delta$ is a Lipschitz manifold, if
  $ \| J - J_0\|_{L^{1}} + \| u_0 - u_{J}\|_{L^{2}} < 2 \delta$.
  Thus, in particular,
    $ U_{J}^{\delta}$ is a   Lipschitz manifold, if
  $ \| J - J_0\|_{L^{1}} + \| \tilde{u}_0 - u_{J} \|_{L^{2}} < {\delta}$,
  for any   $ \tilde{u}_{0} \in E_{J_0} $
 with   $ \|  \tilde{u}_0 - u_{0} \|_{L^{2}} < \delta$. Taking a finite subcovering
  of the covering of $E_{J_0}$ by balls $ B(u_0, \delta(u_0))$,
 with $u_0$ varying in $E_{J_0}$, the first part of the result follows
  with $\delta$ chosen as the minimum of those $\delta(u_0)$.

  Now, if  $ \varepsilon >0$ and  $ u_{0} \in E_{J_0} $, there exists,
    by  Lemma \ref{contunstball},   $\delta = \delta(u_0)$ such that,
 if
  $ \| J - J_0\|_{L^{1}} + \| u_0 - u_{J} \|_{L^{2}} < 2 \delta$, then
   \[ \textrm{dist}( U_{J}^\delta , U_{J_{0}}^\delta     ) +
dist( U_{J{0}}^\delta, U_{J}^\delta   )  < \varepsilon/2.
\] If   $ \tilde{u}_{0} \in E_{J_0} $ is such that
    $ \|  \tilde{u}_0 - u_{0} \|_{L^{2}} < \delta$ and
  $ \| J - J_0\|_{L^{1}} + \| \tilde{u}_0 - u_{J}\|_{L^{2}} <  \delta$
 then, since
 $ \| J - J_0\|_{L^{1}} + \| {u}_0 - u_{J}\|_{L^{2}} <  2 \delta$
  \[
   \begin{array}{l}
  \textrm{dist}( U_{J}^\delta(u_{J}) ,
   U_{J_0}^\delta(\tilde{u}_0)     ) +
   \textrm{dist}( U_{J_0}^\delta(\tilde{u}_0),
   U_{J}^\delta(u_{J}))   \\
    \quad   <
  \textrm{dist}( U_{J}^\delta(u_{J}) ,
   U_{J_0}^\delta(u_0)     ) +
   \textrm{dist}( U_{J_0}^\delta({u}_0),
   U_{J}^\delta(u_{J}))  +
   \textrm{dist}( U_{J_0}^\delta(\tilde{u}_0) ,
   U_{J_0}^\delta(u_0)     ) \\
     \quad +
   \textrm{dist}( U_{J_0}^\delta({u}_0),
    U_{J_0}^\delta(\tilde{u}_0))< \varepsilon.
 \end{array}
 \]

  By the same procedure above of  taking a finite subcovering
  of the covering of $E_{J_0}$ by balls $ B(u_0, \delta(u_0))$,
  and  $\delta$  the minimum of those $\delta(u_0)$, we conclude that
  \[ \textrm{dist}( U_{J}^\delta(u_{J}) ,
   U_{J_0}^\delta(\tilde{u}_0)     ) +
   \textrm{dist}( U_{J_0}^\delta(\tilde{u}_0),
   U_{J}^\delta(u_{J})) < \varepsilon
 \]
 if
 $ \| J - J_0\|_{L^{1}} + \| \tilde{u}_0 - u_{J}\|_{L^{2}} <  \delta$,
 for any   $\tilde{u}_0 \in E_{J_0}$. This proves the result claimed.
 \qed

\subsection{Characterization  of the attractor}

 As a consequence of its gradient structure (see Remark 4.7 of \cite{Silva4}),
 the attractor of  the flow generated by  {\bf (P)}$_{J}$ is given by
 the union of the  unstable set of the set of equilibria. We prove below a more precise characterization.

 As is well known in the literature, an  equation of the form
$$
\dot{x}+Bx=g(x),\label{A.1}
$$
where $B$ is a bounded linear operator on a Banach space $X$ and
$g:X\rightarrow X$ is a $C^{2}$ function, may be rewritten in the
form
\begin{equation}
\dot{x}+Ax=f(x), \label{A.2}
\end{equation}
where $A=B-g'(x_{0})$ and $f(x)=g(x_{0})+r(x)$, with $r$
differentiable and $r(0)=0$.

The following result has been proven in \cite{----R}.
\begin{teo}
Suppose the spectrum $\sigma(A)$ contains $0$ as a simple
eigenvalue, while the remainder of the spectrum has real part
outside some  neighborhood of  zero. Let  $\gamma$ be a curve of
equilibria of the flow generated by (\ref{A.2}), of class $C^{2}$.
 Then there exists a
neighborhood $U$ of $\gamma$ such that, for any $x_{0}\in U$ whose
 positive orbit     is
 precompact and whose
$\omega$-limit set $\omega(x_{0})$ belongs to $\gamma$, there exists
a unique point $y(x_{0})\in \gamma$ with $\omega(x_{0})=y(x_{0})$.
Similarly,  for any $x_{0}\in U$ with
  bounded  negative orbit   and
$\alpha$-limit set $\alpha(x_{0})$   in  $\gamma$, there exists a
unique point $y(x_{0})\in \gamma$ such that
$\alpha(x_{0})=y(x_{0})$.\label{Teorema A.1}
\end{teo}


\begin{prop}
 Assume  the hypotheses (H1)-(H4) hold. Let $E_{J}$ be the set of the
equilibria of $T_{J}(t)$. For $u\in E_{J}$, let $W_{J}^u(u)$ be the
unstable set of $u$. Then
$$
{\cal A}_{J}=\bigcup_{u\in E_{J}}W_{J}^u(u).
$$\label{Lema 7.7}
\end{prop}
\dem From Remark 4.7 of \cite{Silva4}, follows that
$$
{\cal A}_{J}=W_{J}^u(E_{J}).
$$
There exists only a finite number, $\{u_{1},\cdots, u_{k}\}$ of
constant
 equilibria since  they are all
hyperbolic. For each  nonconstant equilibrium $u \in E_{J}$, there
is a  curve $M_{u} \subset E_{J}\subset {\cal A}_{J}$. From Lemma
\ref{Lema isolated} these curves $M_{u}$ are all isolated and, since
${\cal A}_{J}$ is compact, it follows that there exists only a
finite number of them; $M_{1}, \ldots, M_{n}$. Thus
$$
{\cal A}_{J}=\bigg(\bigcup_{i=1}^{n}W_{J}^u(M_{i})\bigg) \bigcup
\bigg( \bigcup_{j=1}^{k}W_{J}^u(u_{j}) \bigg).
$$
From Theorem \ref{Teorema A.1} follows that
$$
W_{J}^u(M_{i})=\bigcup_{v\in M_{i}}W_{J}^u(v), \,\,\, i=1,\cdots, n.
$$
Therefore
$$
{\cal A}_{J}=\bigcup_{v\in E_{J}}W_{J}^u(v),
$$
which concludes the proof. \qed

 \subsection{ Proof of the lower semicontinuity} \label{prooflower}

Using the results obtained in the previous subsections, the proof of
the lower semicontinuity can now be adapted from Lemma 3.8 and
Theorem 3.9 of \cite{Severino2}, as shown below.


\begin{lema}
Assume the same hypotheses of Proposition \ref{Lema 7.7}. Then,
given $\varepsilon
>0$, there exists $T>0$ such that,  for all $u \in {\cal
A}_{J_{0}}\backslash E_{J_{0}}^{\varepsilon}$
$$
T_{J_{0}}(-t)u \in E_{J_{0}}^{\varepsilon},
$$
for some $t\in [0,T]$, where $E_{J_{0}}^{\varepsilon}$ is the
$\varepsilon$-neighborhood of $E_{J_{0}}$. Furthermore, when
$\varepsilon$ is sufficiently small,
$$
T_{J_{0}}(-t)u \in U_{J_{0}}(u_{0}),
$$
for some $u_{0}\in E_{J_{0}}$, where $U_{J_{0}}(u_{0})$ is the local
unstable manifold of $u_{0}\in E_{J_{0}}$. \label{Lema 7.8}
\end{lema}
\dem Let $\varepsilon >0$ be  given and $u\in {\cal
A}_{J_{0}}\backslash E_{J_{0}}^{\varepsilon}$. From Proposition
\ref{Lema 7.7}, it  follows that
$$
u\in W_{J_{0}}^{u}(\bar{u})\backslash E_{J_{0}}^{\varepsilon}.
$$
for some $\bar{u}\in E_{J_{0}}$. Thus, there exists
$t_{u}=t_{u}(\varepsilon)< \infty$ such that $ T_{J_{0}}(-t_{u})u
\in E_{J_{0}}^{\varepsilon}. $ By continuity of the operator
$T_{J_{0}}(-t_{u})$, there exists $\eta_{u}
>0$ such that
$ T_{J_{0}}(-t_{u})B(u, \eta_{u}) \subset E_{J_{0}}^{\varepsilon}, $
where $B(u, \eta_{u})$ is the ball of center  $u$ and radius
$\eta_{u}$. By compactness, there are $u_{1}, \cdots, u_{n} \in
{\cal A}_{0}\backslash E_{J_{0}}^{\varepsilon}$ such that
$$
{\cal A}_{J_{0}}\backslash E_{J_{0}}^{\varepsilon} \subset
\bigcup_{i=1}^{n}B(u_{i},\eta_{u_{i}}),
$$
with $T_{J_{0}}(-t_{u_{i}})B(u_{i}, \eta_{u_{i}}) \subset
E_{J_{0}}^{\varepsilon}$, for $i=1, \ldots, n$. Let
$T=\max\{t_{u_{1}},\cdots, t_{u_{n}}\}$.  Then,  for any  $u\in
{\cal A}_{J_{0}}\backslash E_{J_{0}}^{\varepsilon},$ $
T_{J_{0}}(-t)u\in E_{J_{0}}^{\varepsilon}, $ for some $t\in [0,T]$.
Since $u\in W_{J_{0}}^{u}(\overline{u})\backslash
E_{J_{0}}^{\varepsilon}$, for some $\overline{u}\in E_{J_{0}}$ and
$T_{J_{0}}(-t)u\in E_{J_{0}}^{\varepsilon}$, to  conclude that
$T_{J_{0}}(-t)u\in U_{J_{0}}(\bar{u})$, when $\varepsilon$ is
sufficiently small, it is  enough to show that there exists $\delta
>0$ such that $W_{J_{0}}^{u}(v)\cap B(v,\delta)\subset
U_{J_{0}}(v)$, for all $v\in E_{J_{0}}$. Therefore,  the conclusion
follows immediately from Lemma
 \ref{contunstball}. \qed

\begin{teo}
Assume the hypotheses (H1)-(H4). Then the family of attractors
$\{{\cal A}_{J}\}$ is lower semicontinuous with respect to the
parameter $J$ at $J_{0}\in {\cal J}$.\label{Teorema 7.9}
\end{teo}
\dem Let $\varepsilon >0$ be given. From Lemma \ref{Lema 7.8}, there
is $T>0$ such that,  for all $u\in {\cal A}_{J_{0}}\backslash
E_{J_{0}}^{\varepsilon}$, there exists $t_{u}\in [0,T]$ such that
\begin{equation}
\bar{u}:= T_{J_{0}}(-t_{u})u \in U_{J_{0}}(u_{0}), \label{7.28}
\end{equation}
for some $u_{0}\in E_{J_{0}}$.
 Since $T_{J_{0}}(t)$ is a continuous family of
  bounded operators, there exists $\eta
>0$ such that,  for all  $t\in [0,T]$
\begin{equation}
\|z-w\|_{L^{2}} < \eta \Rightarrow \|T_{J_{0}}(t)z -
T_{J_{0}}(t)w\|_{L^{2}} < \frac{\varepsilon}{2}.\label{7.29}
\end{equation}

By the uniform  continuity of the equilibria and
 local unstable manifolds with respect to
the parameter $J$ asserted by Theorem \ref{Teorema 7.6} and Lemma
 \ref{contunstballunif} , there exists $\delta^{*}>0$ (independent of
 $u$) such that $\|J-J_{0}\|_{L^{1}} < \delta^{*}$ implies the
 existence of $u_{J} \in E_{J}$ and
  some $\bar{\bar{u}}_{J}\in U_{J}(u_{J})$
with
\begin{equation}
\|\bar{\bar{u}}_{J}-\bar{u}\|_{L^{2}} < \eta, \label{7.30}
\end{equation}
where $U_{J}(u_{J})$ denotes the local unstable manifold of the
equilibrium $u_{J}$ of $T_{J}(t)$. Thus, when $\|J-J_{0}\|_{L^{1}} <
\delta^{*}$  we obtain, from (\ref{7.29}) and (\ref{7.30})
\begin{equation}
\|T_{J_{0}}(t)\bar{\bar{u}}_{J}-T_{J_{0}}(t)\bar{u}\|_{L^{2}}
<\frac{\varepsilon}{2}
 \quad \text{for any}  \quad t \in [0,T].\label{7.31}
\end{equation}

  On the  other hand, from continuity of the flow
  with respect to parameter $J$, (see Lemma 10 of \cite{Silva3}),
there exists $\overline{\delta}>0$ such that $\|J -J_{0}\|_{L^{1}} <
\overline{\delta}$ implies
 \begin{equation}
\|T_{J}(t)(u)-
 T_{J_{0}}(t)(u) \|_{L^{2}}
<\frac{\varepsilon}{2},\label{7.32}
\end{equation}
 for any
 $u \in B(0,2 \tau\|J\|_{\infty}S_{max}+h)$) and $t\in [0,T]$. In particular, (\ref{7.32}) holds for
 $u= \bar{\bar{u}}_{J}$ and $t= t_u$.

Choose $\delta=\min\{\delta^{*}, \overline{\delta}\}$ and
 let $v_{J}:=T_{J}(t_{u})\bar{\bar{u}}_{J}$.
  Note that
  $ v_{J} \in {\cal
  A}_{J}$, since  $\bar{\bar{u}}_{J} \in U_{J}(u_{J}).$

 Thus, using (\ref{7.31}) and (\ref{7.32}) we obtain, when $\|J
-J_{0}\|_{L^{1}} < \delta$
\begin{eqnarray*}
\|v_{J}-u\|_{L^{2}}&=&\|T_{J}(t_{u})\bar{\bar{u}}_{J}-T_{J_{0}}(t_{u})
 \bar{u} \|_{L^{2}}\\
&\leq &
\|T_{J}(t_{u})\bar{\bar{u}}_{J}-T_{J_{0}}(t_{u})\bar{\bar{u}}_{J}\|_{L^{2}}
+
\|T_{J_{0}}(t_{u})\bar{\bar{u}}_{J}-T_{J_{0}}(t_{u})\bar{u}\|_{L^{2}}\\
& < & \varepsilon.
\end{eqnarray*}

When $u\in E_{J_{0}}^{\varepsilon}\subset {\cal A}_{J_{0}}$ this
conclusion follows straightforwardly  from the continuity of
equilibria. Thus
 the lower semicontinuity of attractors follows.
\qed

\section{A concrete example}

In this section we illustrate the results of the previous sections
to the particular case of (\ref{WC}) where
 $f(x)=(1+e^{-x})^{-1}$ and
$$
\widetilde{J}(x)=\left\{\begin{array}{ccccc}
e^{\frac{-1}{1-x^{2}}},  \,\,if \,\, |x|< 1,\\
                      0,\,\, if \,\,  |x| \geq 1.\\
\end{array}
\right.
$$
The function $f$ has been motivated by similar functions in 
\cite{Coombes}, \cite{Kubota} and \cite{Wilson} and the function  $\widetilde{J}$ has been adapted from a test function
in \cite{Butkov}.

In this case, we can rewrite equation (\ref{WC}) as
\begin{equation}
\frac{\partial v(x,t)}{\partial t}=-v(x,t)+
\int_{-1}^{1}e^{\frac{-1}{1-(x-y)^{2}}}(1+e^{-v(y)})^{-1}dy+ h
\label{Ex1}.
\end{equation}

As mentioned in the introduction, defining $\varphi:
\mathbb{R}\rightarrow S^{1}$ by
$\varphi(x)=exp^{i\frac{\pi}{\tau}x}$ and, for $v\in
\mathbb{P}_{2\tau}$, $u:S^{1}\rightarrow \mathbb{R}$ by
$u(\varphi(x))=v(x)$ and writing
$J(\varphi(x))=\widetilde{J}^{\tau}(x)$, where
$\widetilde{J}^{\tau}$ denotes the $2\tau$ periodic extension of the
restriction of $\widetilde{J}$ to interval $[-\tau,\tau]$, $\tau
>1$, the equation (\ref{Ex1}) is equivalent to
equation
\begin{equation}
\frac{\partial u(w,t)}{\partial t}=-u(w,t)+
\int_{S^{1}}J(wz^{-1})(1+e^{-u(z)})^{-1}dz+ h, \label{Ex2}
\end{equation}
with now $dz=\frac{\tau}{\pi}d\theta$, where $d\theta$ denotes
integration with respect to arc length.

\subsection{Check hypotheses}

The function $f$  satisfies the hypotheses (H1) and (H2) and (H4),
with $k_{1}=S_{max}=1$, $L=\ln2$ and $ k_{2}=\frac{1}{2}$ in
(\ref{1.4}) and  the function $J$ satisfies the hypothesis (H3)-$b$
assumed in the Section 3.

In fact, note that $f'(x)=(1+e^{-x})^{-2}e^{-x}>0$. Then, since $1<
(1+e^{-x})^{2}\leq 4$, $\forall \,\, x \in \mathbb{R}$, follows that
$$\frac{1}{4}\leq (1+e^{-x})^{-2} <1.$$
Thus
$$
|f(x)-f(y)|< |x-y|.
$$
In particular, since $f(0)=\frac{1}{2}$, we have
$$
|f(x)|< |x|+\frac{1}{2}, \,\, \forall \, x\in \mathbb{R}.
$$
Furthermore, since
$f''(x)=2(1+e^{-x})^{-3}e^{-2x}-(1+e^{-x})^{-2}e^{-x}$, we have
$|f''(x)| < 3$, $\forall \,\, x \in \mathbb{R}$, it implies that
$f'$ is locally Lipschitz. Hence (H1) and (H4) are satisfied.

To verify (H2), we begin by noting that
  $0<|(1+e^{-x})^{-1}|<1$ and
$f^{-1}(x)=-\ln(\frac{1-x}{x})$. Thus by a direct computation we
obtain that, for $0\leq s\leq 1$,
$$
\left|\int_{0}^{s}-\ln(\frac{1-x}{x})dx\right|\leq \ln 2.
$$

Finally, to verify (H3), fix a equilibrium solution $u_{0}$  of
(\ref{1.1}), then from Remark \ref{autovalor}
$$
u_{0}'=J*((f'(u_{0})u_{0}')),
$$
that is, zero is eigenvalue of $DF_{u}(u_{0})$ with eigenfunction
$u'_{0}$. Now, from Remark \ref{Remark 3.1}, $DF_{u}(u_{0})$ is
self-adjoint operator. Then, to prove that zero
 is simple eigenvalue, it is enough to show that if
$v\in Ker(DF_{u}(u_{0}))$ then, $v=\lambda u_{0}$ for some $\lambda
\in \mathbb{R}$.

For this, let $v\in L^{2}(S^{1})$ be such that $DF_{u}(u_{0})(v)=0$.
Then
$$
v=J*((f' \circ u)v).
$$
Hence, using Holder inequality, for any $\lambda \in \mathbb{R}$, we
have
\begin{eqnarray*}
|v(w)-\lambda u_{0}'(w)|&=& |J*[f'(u_{0})v- \lambda f'(u_{0})u_{0}^{'}](w)|\\
&\leq&|J*[f'(u_{0})v- f'(u_{0})\lambda u_{0}^{'}](w)|\\
&\leq&\sqrt{2\tau}\|J\|_{\infty}\|f'(u_{0})v- f'(u_{0})\lambda
u_{0}^{'}\|_{L^{2}}.
\end{eqnarray*}
But
\begin{eqnarray*}
\|f'(u_{0})v-f'(u_{0})\lambda u_{0}')\|_{L^{2}}&=&\|f'(u_{0})[v-\lambda u_{0}']\|_{L^{2}}\\
&<& k_{1}\|v-\lambda
u_{0}^{'}\|_{L^{2}}\\
&=&\|v-\lambda u_{0}^{'}\|_{L^{2}}.
\end{eqnarray*}
Then
\begin{eqnarray*}
|v(w)-\lambda u_{0}'(w)|\leq \sqrt{2\tau}\|J\|_{\infty}\|v-\lambda
u_{0}^{'}\|_{L^{2}}.
\end{eqnarray*}
Now, since $0\leq \widetilde{J}(x) \leq e^{-1}$, follows that
$\|J\|_{\infty}\leq \frac{1}{e}$. Thus
\begin{eqnarray*}
\|v-\lambda u_{0}^{'}\|_{L^{2}}\leq  \frac{2\tau}{e}\|v-\lambda
u_{0}^{'}\|_{L^{2}}.
\end{eqnarray*}
It implies
$$
(1-\frac{2\tau}{e})\|v-\lambda u_{0}^{'}\|_{L^{2}}\leq 0.
$$
Thus, choosing $\tau$ such that $\frac{2\tau}{e}<1$, follows that
$v=\lambda u_{0}^{'}$ in $L^{2}(S^{1})$. Hence, zero is simple
eigenvalue of $DF_{u}(u_{0})$.

Therefore all results of Sections 2 and 3 are valid for the flow
generated by equation (\ref{Ex2}).

\subsection{Concluding remarks}

\begin{Remark}
In (\ref{Ex1}), this choice for $\widetilde{J}$  implies that we are
in the case of lateral-inhibition type fields (short-range
excitation and long-range inhibition),  (see for example,
\cite{Coombes} \cite{Kishimoto}  and \cite{Rubin}). Similar
connection functions ( type "Mexican hat" ) as
$\widetilde{J}(x)=e^{-a|x|}$, $a>0$,
$\widetilde{J}(x)=2\sqrt{\frac{b}{\pi}}e^{-bx^{2}}$, $b>0$ or
$\widetilde{J}(x)=e^{-a|x|}-e^{-b|x|}$, $0< a<b$,  has been used
often  in previous work, (see, for example, \cite{Ermentrout},
\cite{Ermentrout2}, \cite{Pinto}, \cite{Pinto2} and \cite{Rubin}).
Hoping to make the model more realistically the connectivity
existing in the prefrontal cortex, in \cite{Laing} is considered the
synaptic connection function
$\widetilde{J}(x)=e^{-b|x|}(b\sin|x|+\cos x))$, which changes sign
infinitely often.
\end{Remark}

\begin{Remark}
Note that, the equivalence between the equations (4.1) and (4.2), given in the formulation above, implies that  the  lateral-inhibition type connectivity function (short-range excitation and long-range inhibition) in (\ref{Ex1}), when restrict to space of $2\tau$-periodic functions, results in a recurrent-excitation type connectivity function in (\ref{Ex2}).  Therefore,  thus as in \cite{Laing}, we hope have a connectivity function $J$ that represents more realistically the connectivity existing in brain
activities, since it is  known that electrical discharges from brain cells result in a recurrent seizure disorder such as migraine and epilepsy (see, for example, \cite{Ruktamatakul}).
\end{Remark}


\section*{Acknowledgements}

The author would like to thank the anonymous referee for his/her
reading of the manuscript and valuable suggestions. Our gratitude
goes also to Oxford International English that checked  the English
usage this paper.


\end{document}